\documentclass[11pt]{article}
\usepackage[utf8]{inputenc}
\usepackage{graphicx}
\usepackage{amsfonts}
\usepackage[width=17cm, left=1cm]{geometry}
\usepackage{amssymb}
\usepackage{amsmath}
\usepackage{amsthm}
\usepackage[russian]{babel}

\newtheorem{definition}{Definition}
\newtheorem{proposition}{Proposition}
\newtheorem{theorem}{Theorem}

\newtheorem{lemma}{Lemma}
\newtheorem{corollary}{Corollary}

\title{\textbf{COMBINATORIAL ENCODING OF CONTINIOUS DYNAMICS, AND TRANSFER ON THE SPACE OF PATHS} \thanks{Partially supported by RSСF 17-71-20153.}}
\bigskip

\begin{document}

\bigskip

\author{A.~М.~Vershik\thanks{Saint-Petersburg branch of Stekloff Mathematical Institute of Russian Academy of Science, Saint-Petersburg State University, Mathematical Department,Moscow Institute of the Problems of Transmission of information.}}
\date{01.04.2019}
\maketitle
\rightline{TO THE MEMORY OF SERGIY KOLYADA}

\bigskip

\bigskip

\centerline{Abstract}
These notes follow my articles \cite{V,V1}, and give some new important details. We propose here a new combinatorial method of encoding of measure spaces with measure preserving transformations, (or groups of transformations) in order to give new, mostly locally finite geometrical models for investigation of dynamical properties of these objects.

\medskip

\section{Partitions of the infinite-dimensional cube into Weyl simplices and distinguishability problem}

We define an increasing (it means reducing) invariant sequence of finite cylindric measurable partitions
$\{\eta_n\}_n$
of infinite dimensional cube $(I^{\infty}$ with product measure $m^{\infty})$ with Lebesgue measure $m$ on the interval $[0;1]$ as a multiplier. More exactly, the partition $\eta_n$ is the \emph{cylindric partition} which has as a base the partition of the $n$-dimensional cube $I^n$ into open \emph{Weyl simplices}. The meaning of the last sentence is the following:
we consider the space ${\Bbb R}^n$ as a Cartan subalgebra of Lie algebra of series $(A_n)$ with fixed root system. We call "open Weyl simplices"
the intersections of the open Weyl cameras in usual sense with the unit cube $[0,1]^n$ . We also suppose that the correspondence between Weyl cameras and elements of Weyl group (which is symmetric group $S_n$) is also fixed. Although the language of the theory of Lie algebras it is not necessary for description in our example, we nevertheless use it in order to generalise in future the construction to the case of simple Lie algebras of other series.

 We will consider the subset of the cube $I^n$ of full (Lebesgue) measure $m^n$  that consists of the vectors with pairwise distinct coordinates.

Each Weyl simplex is identified with some permutation of $\textbf{n} = \{1; 2; \dots n\}$ -element of Weyl group.

Recall that the \emph{finiteness} of partitions $\eta_n$ means that number of elements of positive measure is finite; in our
case it equals to $\mbox{ord}(S_n) = n!$.

The \emph{increase} of the sequence $\{\eta_n\}_n$ means that for all $n$ each element $C$ of partition
$\eta_{n+1}$
is (as a subset of $I^n$) a subset of element $D$ of the partition $\eta_n$, and each element $D \in \eta_n$
(as subset of $I^n$) is exactly a disjoint union of
all those elements of partition $\eta_{n+1}$ which belong to it. The union of all elements of all partitions $\eta_n, n=1,2 \dots$ (the set of all Weyl simplices) generates \emph{a tree which we denote  $W$} and call the permutation tree (or factorial tree).

The $n$-th level of the tree $W$ has $n!$ vertices  which correspond to the elements of partition $\eta_n$, and each vertex can be identified with a permutation from $S_n$.

The edges of the tree join each element  $C \in \eta_{n+1}$ with the element $D$ of partition $\eta_n$ that contains element $C\in D$; see the picture 1.
Let us define the coordinates for Weyl simplecies and the corresponding permutation as follow.
Let $x^n=(x_1,x_2,\dots x_n)$ be a point of $I^n$, then the coordinates

$$(k_1,k_2,\dots, k_n)$$ corresponding to the simplex $\sigma_{x^n}$
are given by the formula
    $$k_i=\#\{s \in {\bf n}: x_s<x_i\}, \quad i=1,2, \dots n.$$
In particular coordinates $(k_1,k_2,\dots, k_n)$ do not depend on the choice of a point of the simplex.

\begin{figure}
\includegraphics[width=\linewidth]{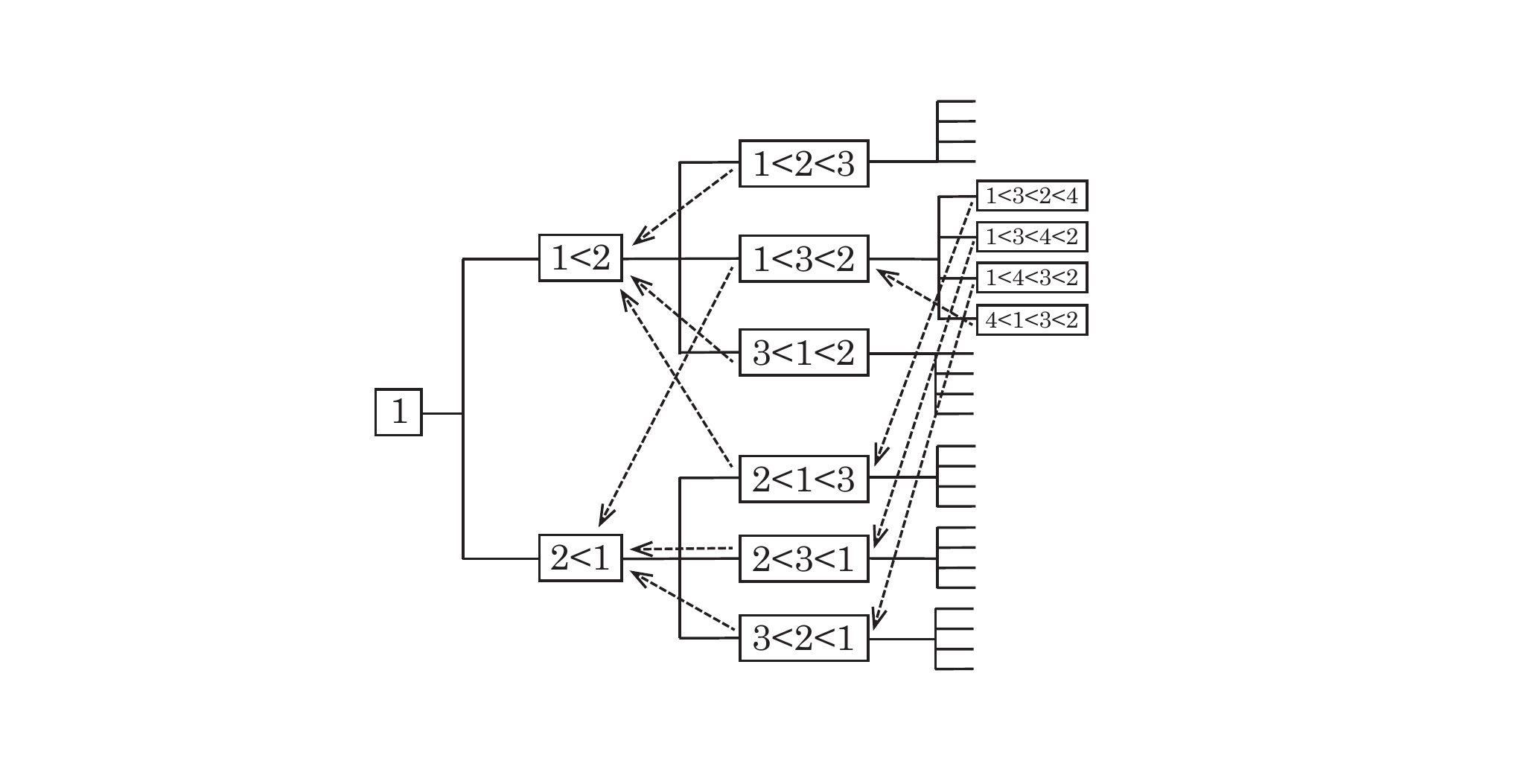}
\caption{Tree of the permutations with translation.} \label{fig:Tree}
\end{figure}

The \emph{invariance} of the sequence of partitions $\{\eta_n\}_n$ with respect to one-sided shift $S$ of the cube $I^{\infty}$ means that the images
$Sx, Sx'$ of almost all pairs points $x, x'$ of $I^{\infty}$, which belong to the same element of partition $\eta_n$ (or belong to the same Weyl simplex)
will belong to the same element of partition (or belong to the same Weyl simplex) of $\eta_{n-1}$, for $n = 2,3 \dots$.

This means that we constructed another system of edges on the tree $W$: each element of the partition $\eta_n$ that corresponds
to the unique vertex of level $n$ has an edge which goes to the vertex of the previous level $n-1$ that corresponds to some element of partition $\eta_{n-1}$. We call \emph{translations} such edges between permutations.

 \begin{proposition}
The permutation $(r_1,r_2,\dots,r_{n-1})$ corresponding to the shifted simplex
$\Sigma_{Sx^n}$ is given by the formula
      $$
 r_i=\begin{cases}
 k_{i+1} &\text{if $k_{i+1}<k_1$},\\
 k_{i+1}-1&\text{if $k_{i+1}>k_1$}.
 \end{cases}
 $$
 \end{proposition}

 The proof immediately follows from definitions.
The step from the sequence $\{k_i\}$ to the sequence $\{r_i\}$ is the "dynamics"\ of our coordinates.

So we equip our sequence of partitions $\{\eta_n\}_n$ with the tree with additional structure ---translation and transfer (see picture 1). The tree $W$ with these structures we call the \emph{skeleton} of the permutation tree.

 An ordinary edge of the tree joins the permutation $g \in S_n$ with a permutation
$h \in S_{n+1}$ if and only if $g$ is a result of the removal of $(n+1)$-st object from the permutation $h$. In opposite to this the translation joins the permutation $g \in S_n$
with permutation $f \in S_{n-1}$ if $f$ is a result of the removal of the first object of permutation $g$
\footnote{Here abusing notation we have identified a permutation as an element of the symmetric group $S_n$ and a permutation which is the image of the natural order $(1,2, \dots n)$ under this permutation; of course one must distinguish these notions, see \cite{V1}}

Thus we have the skeleton which is a combinatorial scheme of the increasing shift-invariant sequence of
finite partitions. Now we want to consider this data as the new coordinates of the points of the cube.

\begin{definition}(Transfer)

Consider the space $T(W)$ of all paths of the permutation tree $W$. The operation $\lambda$ that sends the permutation $g\in S_n$ of level $n$ by translation to permutation $h\in S_{n-1}$ of level $n-1$ induces the map $\Lambda:T(W)\rightarrow T(W)$: the path $\{g_1,g_2,g_3 \dots\}$ goes to the path $\{\lambda(g_2),\lambda(g_3), \dots \}$; we call this operation $\Lambda$ \emph{transfer} on the tree $W$. Because $g\prec h$ gives $\lambda(g)\prec\lambda(h)$   (two removals commute), this operation defines the mapping of the paths of our tree.
\end{definition}

We will give a transparent formula for this transfer,
see Proposition 2 below.

For $x=(x_1,x_2,x_3, \dots \in I^{\infty} $
let us consider the quantities:
$t_n(x)=\sharp\{k\leq n: x_k<x_n\}$
and the map
$$ x\mapsto t(x)=\{t_1(x),t_2(x),\dots \},\quad t_n(x)\in \textbf{n},\quad n=1,2 \dots$$

Let us correspond to a given point $x=\{x_1,x_2, \dots x_n \dots\}\in I^{\infty}$ and to corresponded Weyl simplex the permutation as follow sequence of the elements of the partitions $\eta_n$  which contains point $x$:
because each element of $\eta_n$ corresponds to a permutation $g \in S_n$ we have sequence
$$ x\mapsto \{t_1(x),t_2(x) \dots t_n(x) \dots \}.$$

It is clear that this sequence $\{t_n(x)\}$ is a path in the skeleton of our tree $W$.

This is the combinatorial encoding of the infinite dimensional unit cube $I^{\infty}$ with respect to the partition into the Weyl simplices.

\begin{definition}
Denote compact space:
$$\mathfrak{M}=\{\{t_n\}_{n=1}^{\infty}, \quad t_n\in \textbf{n}, \quad n=1,2 \dots\} $$
and call $\mathfrak{M}$ the \emph{triangular compact set}.
\end{definition}

We define the map
  $J:I^{\infty}\rightarrow \mathfrak{M}; x=\{x_1,x_2, \dots\}\mapsto \{t_1(x),t_2(x) \dots\}:$

 $$J:[0,1]^{\infty}\rightarrow \mathfrak{M},\quad  x=\{x_n\}\in I^{\infty}; \quad
    J(\{x_n\})=\{t_n=t_n(x_1,\dots, x_n)\}_{n=1}^{\infty}; \quad\quad t_n=\sharp\{i:1\leq i \leq n \quad x_i<x_n\}.$$

  Let us define the probability measure  $\mu$ on the compact set $\mathfrak{M}$ that is the product measure
  of the uniform measures on each multiplier $\textbf{n}$.

\section{Distinguishability of Weyl encoding}

  Our first observation is the following:

  \begin{theorem}
  The map $J$ is the isomorphism between measure spaces
  $(I^{\infty}, m^{\infty})$ and $(\mathfrak{M},\mu);$
  this means that coordinates $\{t_n(x)\}_{n=1}^{\infty}$ define almost all $x$, in other terms: the product of partitions
  $\eta_n$ is the \emph{identity} partition or the partition of the infinite dimensional cube $I^{\infty}$ into the separate points; we denote it by $\epsilon$:

     $$\bigvee_{n=1}^{\infty}\eta_n=\epsilon.$$
  \end{theorem}

 In geometrical terms this means that the set of all Weyl simplices of the cube $I^{\infty}$ separates almost all points of the cube. More paradoxical formulation is the following: almost every
 sequence $\{x_1,x_2, \dots\}$ of the cube $I^{\infty}$ with respect to Lebesgue measure $m^{\infty}$ can be restored if someone knows only all inequalities $x_n>x_m$ or $x_n<x_m$ for all $n,m\in \Bbb N$. Or even in more expressive form: \emph{almost every Weyl simplex consists of the unique point, and almost every Weyl camera consists of the unique ray}.

 \begin{proof}
 Let two sequences $\{x_n\}$ and $\{x'_n\}$ have the same inequalities for all pairs
 $x_n>x_m \Leftrightarrow x'_n>x'_m$.
 Suppose that $x_1\ne x'_1$, then there exists $k$
 such that $x_k \in (x_1,x'_1), x'_k\in (x_1,x'_1)$,
 so the pair $x_1$ and $x_k$ has the opposite inequality than  $x'_1$ and $x'_k$ --- a contradiction.
 \end{proof}
 More precise from of this claim is the following.

\begin{lemma}
The limiting partition $\eta=\lim_n \eta_n$  of the infinite-dimensional cube $I^{\infty}$ (the limit of the partitions into open Weyl simplices) coincides
 $\bmod\, 0$ (with respect to the Lebesgue measure) with the partition into singletons. In other words, the distinguishability problem for the partition into Weyl simplices has a positive answer. Therefore, the map $J$ is an isomorphism of measure spaces. In more detail, there exists a set of full Lebesgue measure in $I^{\infty}$ such that for any two points~$\{\xi_n\}$ and~$\{\xi'_n\}$ of this set there exist indices
$i$ and $j$ for which the corresponding coordinates satisfy the opposite inequalities:
$$\xi_i > \xi_j, \quad \mbox{but}\quad \xi'_i <\xi'_j.$$
\end{lemma}

In a somewhat paradoxical form, the lemma can be stated as follows: almost every (with respect to the Lebesgue measure) infinite sequence of points from the interval $[0,1]$ can be uniquely recovered from the list of pairwise inequalities between its coordinates.
Or, even more paradoxically: almost every infinite-dimensional Weyl simplex consists of a single point.\footnote{If, instead of the cubes $I^n,I^{\infty}$  we consider the spaces ${\Bbb R}^n, {\Bbb R}^{\infty}$ with the standard infinite-dimensional Gaussian measures, then our statement looks as follows: almost every infinite-dimensional Weyl chamber consists of a single ray.}

We can say that sequence of the partitions on the Weyl simplices give the positive solution of the \emph{Distinguishability Problem}, --- the problem how to distinguish the points of the infinite dimensional unit cube with the help of some sequence of the finite measurable partitions. The real explanation of such  effect lies in the individual ergodic theorem or in the theorem about uniform distribution on the interval of almost all points of the cube. This mean that the  Lebesgue measure $m^{\infty}$ could be changed to any measure with this property.

\smallskip

Now we must describe the translation and transfer which was defined above in terms of coordinates $\{t_n\}$. This means the rule of changing of these coordinates when the point $x=\{x_1,x_2,x_3, \dots\}$ changes onto  point $Sx=\{x_2,x_3, \dots\}$. Exactly, we want to express $\{t_1(Sx)=t'_1, t_2(Sx)=t'_2, \dots \}$ in terms of $ \{t_1(x),t_2(x),t_3(x) \dots \}$. Such formula will give the formula for the transformation:

$$\Lambda=JSJ^{-1}: \mathfrak{M} \rightarrow \mathfrak{M},$$ of the triangular compact set $\mathfrak M$.

 This map $\Lambda$ is the transfer of the permutation tree and of triangular compact set as the set of paths of that tree. But we will give a precise formula for $\Lambda$.

  Denote $\xi^n=\{\xi_i\}_{i=1}^n$ and let  $d_n(\xi^n)$ be the number of coordinates in the vector $\xi^n$ that are less than $\xi_1$. It is clear that $d_{n+1}$ is either equal to $d_n+1$ if $\xi_{n+1}<\xi_1$, or equal to $d_n$
  if $\xi_{d+1}>\xi_1$. Let us say that
 some positions in the vector $(t_1,t_2, \dots t_n)$
 are \emph{special}: the first one $t_1=1$ is special, then if number of special positions among first $n$ is $d_n(t)$, then $t_{n+1}$ is special iff
 $t_{n+1}\leq d_n$, n=1,2 \dots.

\begin{proposition}(Formula for $\Lambda= JSJ^{-1}$)

 $\Lambda(\{t_n\})=\{t'_n\}$,
 where
 $$
 t'_n=\begin{cases}
 t_{n+1} &\text{if $t_{n+1}$ is special}\\
 t_{n+1}-1&\text{if $t_{n+1}$ is not special.}
 \end{cases}
 $$
\end{proposition}

     The following formula takes place:

  $$t_{n+1}-t'_n =1-(d_{n+1}(t)-d_n(t)).$$

The proof follows automatically from the previous formulas.

Now consider formula for the inverse map.  From the formula for $d_n$ we have
directly the following:

\begin{theorem}
For almost every trajectory  $\{\xi_n\}_n\in I^{\infty}$ with respect to the measure $m^{\infty}$ on $I^{\infty}$ there exists the limit:

$$\lim_n \frac{d_n}{n}=\xi_1.$$
 \end{theorem}

So we can find the first coordinate from the infinite vector $\{t_n\}$. In the same manner we can find the other coordinates $\{x_2,x_3, \dots x_n\}$ using the transfer.

 This means that we can restore the shift $S$
as a map on $I^{\infty}$ with the help of $\Lambda$.

We establish the isomorphism between two triples
  $$(I^{\infty},m^{\infty}, S)\quad \mbox{and}\quad (\mathfrak{S},\mu,\Lambda).$$
  In particular this means that $\Lambda$ is a Bernoulli automorphism of
 $(\mathfrak{S},\mu)$ in the sense of ergodic theory.

By definition, the translation (see Section~2.3) associates with every vertex of level~$n$ (for $n>1$) a vertex of level~$n-1$ following the rule according to which the simplex
$\sigma_{x^n}$ of sequences starting from a vector
$x^n=(x_1,x_2, \dots, x_n)$  changes when we remove the first coordinate $x_1$ after the application of the shift $S$, that is, pass to the vector
$Sx^n=(x_2,x_3, \dots, x_n)$. Recall that the permutation
$(k_1,k_2,\dots, k_n)$ corresponding to the simplex $\sigma_{x^n}$
is given by the formula
    $$k_i=\#\{s \in {\bf n}: x_s<x_i\}\quad i=1,2, \dots n.$$

 Using this rule, we construct the frame corresponding to the tree of Weyl simplices, see Fig.~1.

 Thus we have defined a translation that is a~map from the set of permutations of length~$n$ to the set of permutations of length~$n-1$. In the next section, where we compute the transfer for this graph, we use this map and interpret it in a~slightly different way.

\section{More complicated case: encoding via RSK-correspondence}

In the previous section we considered the partition $\eta_n$ of the cube $I^n$ into Weyl simplices and constructed the coding of the dynamical system $(I^{\infty}, m^{\infty}, S)$ using this kind of partitions. We obtained the combinatorial version of Bernoulli shift as a transformation on the space of paths of the tree (tree of permutations).

Looking on more complicated examples, suppose that the sequence of partitions $\{\alpha_n\}$ of the cube $I^n$ are more coarse than partitions into Weyl simplices $\eta_n$, namely, each element of $\alpha_n \prec \eta_n$ consists of several Weyl simplices. In this case it possible to loose the property of distinguishability: the limit $\bigvee_n \alpha_n$ could be different from the identity partition. How to refine the problem of distinguishability for such sequences of partitions; how to pick up the case when $\lim_n \alpha_n=\epsilon$?

We will keep in mind as one of important further examples the situation with so called RSK-correspondence. The related question was considered in the paper \cite{KV} in 80-th and recently in the important papers \cite{RS,Sn}.

In the framework of this paper the positive answer to the problem of distinguishability for a combinatorial encoding is equivalent to the fact that the homomorphism of Bernoulli shift to the space of infinite Young tableaux with Plancherel measure and Schutzenberger shift is an isomorphism. The question whether this homomorphism is an isomorphism appeared as a result of the paper \cite{KV} in which a generalization of RSK-correspondence for infinite case was defined; and this question was recently answered in the papers \cite{RS,Sn}. Here we mention another approach to the problem of isomorphism in the spirit of this paper and the article \cite{V} as well as the theory of filtrations \cite{UMN}.

We will use the identification of the Weyl simplices with permutations,
more exactly: permutation $g=(i_1,i_2,\dots i_n)\in S_n$ parameterises the Weyl simplex $\sigma_g$ whose elements $x=\{(x_1,x_2,\dots x_n)\in I^n$ has the same ordering of  its coordinates as the ordering of the coordinates of permutation $g$.

The main property of the finite RSK-correspondence (see \cite{St}) is the set-theoretical isomorphism
between symmetric group $S_n$ and set of all pairs of Young tableaux with the same Young diagrams:

$$ S_n=\coprod_{\lambda \in {\hat S}_n} T^P_{\lambda}\times T^Q_{\lambda},$$

where $T^P,T^Q$ is the set of all Young tableaux with the diagram $\lambda$; indexes $P,Q$ means that a tableau
is either insertion (P), or recording (Q).

The index of any simplex is a
permutation $g\in S_n$, but we consider $g$ (as was written above) as a pair of two Young tableaux of the same diagram  $\lambda_n$ : $g=(t_g^P,t_g^Q$),  (shortly we can say that $g$ is associated to the diagram $\lambda_n$).

Now denote by $\Sigma_{\lambda_n}$ a cylindric set: the union of all simplices $\sigma_g$ with the diagram ${\lambda_n}$.

Define two partitions $\eta_n^{\lambda}$ and $\phi_n^{\lambda}$
of $\Sigma_{\lambda_n}$. The first partition $\eta_n^{\lambda}$ is finite
and cylindric,with the elements of type

$$C_{t_Q}=\bigcup_{t_P}\sigma_{t_P,t_Q},\quad  t_P,t_Q\vdash \lambda.$$

The second partition $\phi_n^{\lambda}$ is
cofinite (all elements are finite sets) and each element of partition is the finite set of points  $\{x_n\}\in I^{\infty}$ with equal coordinates $x_i$ when $i>n$,
and vector $(x_1,x_2,\dots x_n)$ run over all normalized permutations with given tableau $t_P$:

$$D_{t_P}=\bigcup_{g=(t_P,t_Q)}\{x=(x_1,x_2,\dots x_n): x_i=\frac{g(i)}{n}, i=1,2 \dots n\}.$$

\begin{lemma}
The partitions  $\eta_n^{\lambda}$  and $\phi_n^{\lambda}$  of the set $\Sigma_{\lambda_n}$ are mutually independent complements: each two elements of the first and the second partitions intersect in just one point, and the properties of independency for conditional measures also hold.
The partition $\eta_n^{\lambda}$ increases in $n$ and
the partition $\phi_n^{\lambda}$ decreases in $n$.
\end{lemma}

So the parameter of the elements of the partition $\eta_n^{\lambda}$ (partition $\phi_n^{\lambda}$) is the set of  permutations with a given $P$-tableau (correspondingly with a given $Q$-tableau). Recall that in the previous section for partition $\eta_n$ the parameter was one permutation.

Now we can consider the partitions $\theta_n$ ($\theta_n^{\bot}$) of the cube $I^{\infty}$ whose restriction to each set $\Sigma_{\lambda_n}$ is $\eta_n^{\lambda}$ (correspondingly $\phi_n^{\lambda}$).  It is clear that the first
 partition is finite and measurable (number of elements is equal to the number of Young tableaux with the given number of cells), the second is cofinite (number of points in one element is the same number). Evidently the partition $\theta_n$ is coarser than partition $\eta_n$: $\theta_n \prec \eta_n$, and partition $\theta_n^{\bot}$ is finer than orbit partition of symmetric group $S_n$ for all $n$. Moreover
\medskip
$\theta_n \prec \theta_{n+1}, n=1,2, \dots$.
$\theta_n^{\bot}\succ \theta_{n+1}^{\bot}, n=1,2, \dots$

It is useful to give an interpretation of these partitionss
in terms of Knuth equivalence and dual Knuth equivalence, see \cite{St} and Fomin's addendum to the Russian translation of this book.

\bigskip
Question: what is the limit of sequence of increasing partitions $\theta_n$?

 \bigskip
 \bigskip

 The following theorem is equivalent to the remarkable theorem  by D.Romik and P.Sniady \cite{RS,Sn} in which they proved that a homomorphism defined by \cite{KV} is a true isomorphism. That homomorphism in the case of Plancherel measure was defined as a map from $I^{\infty}$ to the space of infinite standard Young tableaux. That homomorphism sends Lebesgue measure to Plancherel measure. The question whether this homomorphism is indeed an isomorphism was open for many years and was solved in \cite{RS,Sn}. We give the formulation of the theorem in terms that are related
 to our approach.

 \bigskip

\begin{theorem}
 The limit of the sequence of the partition $\theta$ on the space $I^{\infty}$
  is the identity partition:
 $$\lim_n \theta_n=\epsilon.$$

 Thus the distinguishability problem (see above) has positive solution.
\end{theorem}

\medskip
\begin{corollary}
The limit of the decreasing sequence of partitions $\theta_n^{\bot}$ is the trivial partition:
$$\lim_n \theta_n^{\bot}=\nu$$
\end{corollary}

\smallskip

Remark that in general the fact that the limit of decreasing sequence $\{\alpha_n\}$ of the partitions is trivial, is only necessary but not sufficient for the conclusion that the limit of the increasing sequence $\{\alpha_n^{\bot}\}$ of the partitions  which are independent complements to the $\{\alpha_n\}$ is equal to the identity partition.

In terms of dual Knuth equivalence relation on the cube $I^{\infty}$ the assertion of the theorem means that
this is $\mod 0$ an identity equivalence relation.

But at the same time the corollary of theorem above is the following:

 \begin{proposition} The direct Knuth equivalence relation (see \cite{St}) on the space $I^{\infty}$ is ergodic; this means that there is no non-constant measurable functions that are constant on the classes of dual Knuth equivalence.
\end{proposition}

\medskip
Our numerical simulation have shown very slow convergence in this situation.

 The approach of the papers \cite{RS,Sn} is based on analysis of limit shape of Young diagrams and so called Schitzenbeger's shift ("jeu de taquin"). More precisely, each infinite Young tableau is a path in the Young graph; Schitzenbeger's shift is a shift or, in our terminology, a transfer of the path (see next section). It is  natural to choose the maximal strictly monotonic subset of the path which is called "nerve" of tableau. The main observation of authors \cite{RS,Sn} is that this nerve as a piece-wise line on the lattice ${\Bbb Z}^2_+$ has a limit at infinity which could be identified  (after normalization) with the point of the limit shape of normalized Young diagram. This is the key step to construct the inverse isomorphism in \cite{RS,Sn}.

\medskip
Our approach is different and based on another limit shape theorem. We will give the details in other place. Namely our proof uses the detail analysis of the behavior of $P$-tableaux. These $P$-tableaux are not semi-standard because the entries are real but not natural numbers. In opposite to $Q$-tableaux there is no strong limit (when $n$ tends to infinity) of $P$-tableaux, and their weak limit is zero tableau. Nevertheless after right normalization there is a limit shape; and existence of it
can be extracted from old results about limit shape of standard Young tableaux (not diagram). The inversion formula in \cite{RS} is based on behaviour of $Q$-tableaux.
Our approach to inversion formula is based on the fact of stabilization of $P$-tableaux for a given $Q$-tableau after normalization. Roughly speaking, there is the isomorphism between the realization of Bernoulli scheme ($I^{\infty}$),
infinite Young tableaux ($Q$-tableaux with Plancherel measure) and
normalized $P$-tableaux.

\bigskip

\section{Definition of a transfer for a graded graph}

Let us define the notion of a "transfer" --- an operation on the space of paths of the graded graph (or multigraph). The main idea is to apply the theory of graded graphs to ergodic problems under consideration.

\medskip
Consider an arbitrary Bratteli diagram, i.e., an $\Bbb N$-graded locally finite graph (or even a multigraph)~$\Gamma$. An infinite tree is an example of such a graph. A path in~$\Gamma$ is an infinite maximal sequence of edges (not vertices!) in which the beginning of each edge coincides with the end of the previous edge. Denote the space of all paths by  $T(\Gamma)$; this is a Cantor-like compactum in the inverse limit topology.

\begin{definition}
 A~transformation
  $$\Lambda: T(\Gamma)\rightarrow T(\Gamma)$$
  is called a \emph{transfer} if it is continuous
   in compact topology of the space of paths, and satisfies the following locality property: for each path $t=\{t_n\}_n\in T(\Gamma)$  its image is the path $\Lambda(t)\equiv\{t'_n\}$ that is defined inductively:
  this means that the edge $t'_n$ depends on the fragment of initial path $\{t_k\}_1^{n+1}$ and on the end of the previous edge $t'_{n-1}$.

  We say that a transfer is Markov if the fragment above reduces to the pair of edges $\{t_k\}_n^{n+1}$.
 \end{definition}

For stationary graphs, in which the sets of vertices of every level (except the first one) are isomorphic and these isomorphisms are fixed, the translation rule depends on nothing: an edge connecting vertices~$a$ and $b$ of levels~$n+1$ and~$n+2$ goes to the edge connecting the vertices~$a'$ and~$b'$ of levels~$n$ and~$n+1$ identified with the vertices~$a$ and~$b$, respectively. In this case the transfer is an ordinary shift.

For trees regarded as graded graphs, the definition of transfer coincides with the definition from section 1.

\begin{definition}
A graded graph with a ``transfer'' operation is called a quasi-stationary graph. The path space of a~quasi-stationary graph, regarded as a topological Markov compactum, will be called a quasi-stationary Markov compactum.
\end{definition}

Thus we have described a new type of realizations of automorphisms and endomorphisms with infinite entropy as transfers on quasi-stationary Markov compacta, i.e.
in the space which is locally finite.

According to our definition, a transfer is a shift of sequences of edges, and not of sequences of vertices as in the stationary case. Thus this notion opens new possibilities for realizations of transformations.

\begin{figure}
  \includegraphics[width=\linewidth]{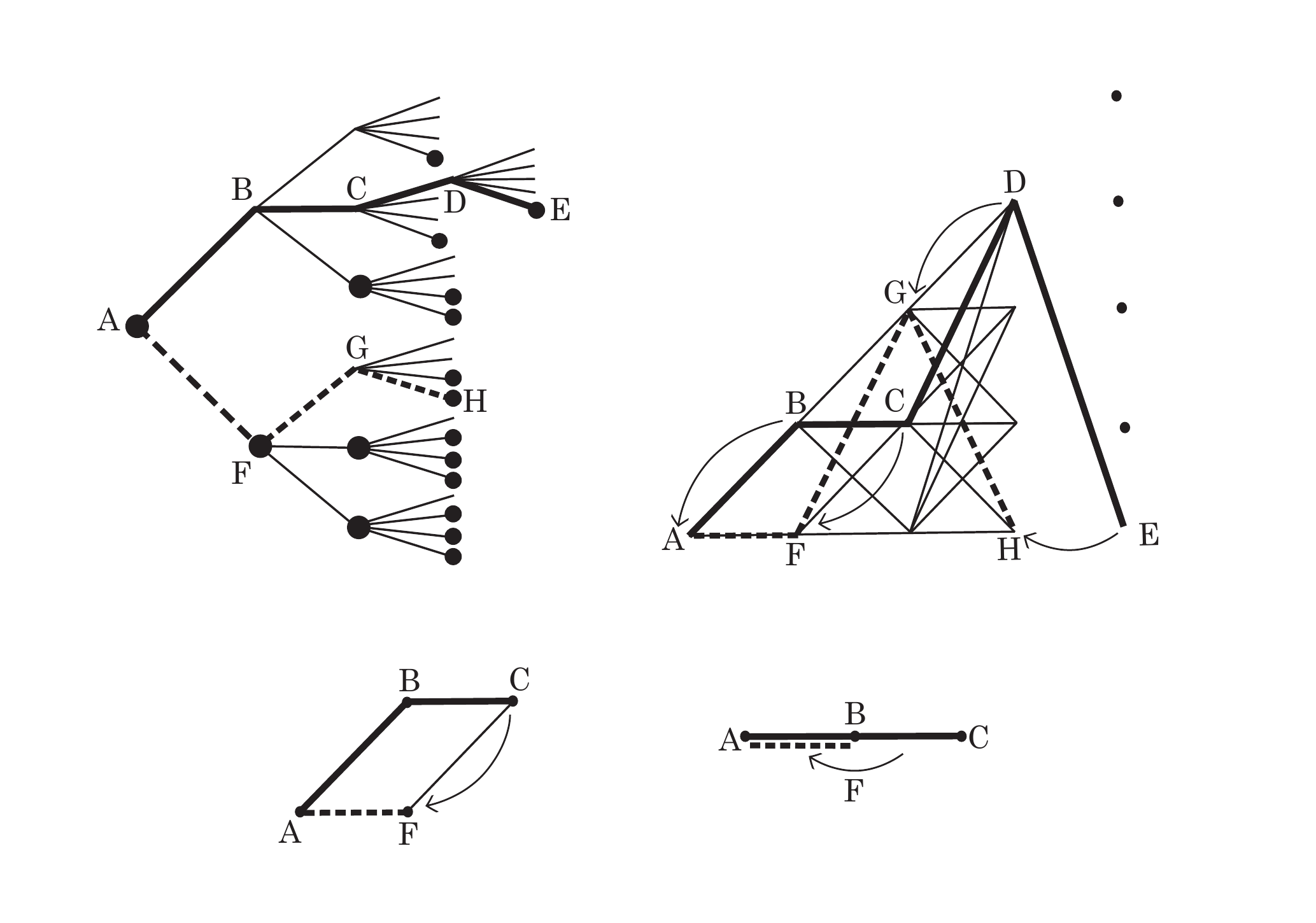}
  \caption{The dashed path is the transfer of the bold one.}
  \label{fig:transfer}
\end{figure}

A transfer defines an additional structure on the graph and, in general, is not uniquely determined by the graph itself, although in some cases there exists a distinguished transfer. For example, the extremely important graphs are the graphs in which every $2$-interval contains either one or two intermediate vertices (this is true for Hasse diagrams of arbitrary distributive lattices, in particular, for Young graph). In this case the natural definition of transfer is Markov transfer, see picture $2$ from which our definition should be clear.

\begin{proposition}
For graphs which are the Hasse diagrams of a distributive lattice of finite ideals of a locally finite countable partial ordered sets with the minimal element, there is a distinguished Markov transfer.  In the case of lattice ${\Bbb Z}_2$ as a poset and Young graph as a graded graph, the transfer coincides with the well-known Sch\"utzenberger transformation,("jeu de taquin" see \cite{St}), so our definition is a generalization of "jeu de taquin" for distributive lattices.
\end{proposition}

The proof follows from a detailed analysis of the definition of transfer.

If  a transfer is defined on the path space of a graded graph, then this space should be regarded as a~nonstationary (or {\it quasi-stationary}) Markov chain, meaning that  the transfer is an analog of the shift. If we have a central measure on the path space that is invariant under the transfer, then we obtain a~quasi-stationary Markov chain with an invariant measure.  Hence the theory of transfer becomes part of ergodic theory, as a nonconventional realization of measure-preserving transformations. We return to all this facts elsewhere.
 
 The author is grateful to Dr.P.Nikitin for preparing of the pictures and for the design of the manuscript.

 \end{document}